\documentclass[10pt,reqno]{article}
\usepackage{amssymb,amsmath,amsthm}
\usepackage[pdftex]{graphicx}
\usepackage{hyperref}

\title{On the quasisymmetric H\"older-equivalence problem for Carnot groups}
\author{P. Pansu\footnote{P. P.~is supported by MAnET Marie Curie Initial Training Network and Agence Nationale de la Recherche grant ANR-2010-BLAN-116-01 GGAA.
}}

\newtheorem*{stheorem*}{Theorem}
\newtheorem{thm}{Theorem}
\newtheorem{prop}[thm]{Proposition}
\newtheorem{lem}[thm]{Lemma}

\newtheorem{cor}[thm]{Corollary}
\newtheorem{defi}[thm]{Definition}
\newtheorem{rem}[thm]{Remark}
\newtheorem{exa}[thm]{Example}
\newtheorem{que}{Question}
\newenvironment{pf}{\begin{trivlist}\item[]{\bf Proof\ }}
{\mbox{}\hfill\rule{.08in}{.08in}\end{trivlist}}

\def\R{\mathbb{R}}
\def\N{\mathbb{N}}
\def\C{\mathbb{C}}
\def\H{\mathbb{H}}
\def\O{\mathbb{O}}

\begin{document}
\maketitle

\abstract{A variant of Gromov's H\"older-equivalence problem, motivated by a pinching problem in Riemannian geometry, is discussed. A partial result is given. The main tool is a general coarea inequality satisfied by packing energies of maps.}

\setcounter{tocdepth}{1}
\tableofcontents

\section{The problem}

\subsection{The H\"older equivalence problem}

In \cite{Gromovcarnot}, M. Gromov conjectured that if there exists a $C^{\alpha}$ homeomorphism from an open set in Euclidean space $\R^3$ to Heisenberg group equipped its Carnot-Carath\'eodory metric, then $\alpha\leq \frac{1}{2}$. This is still open. Gromov proved the upper bound $\alpha\leq \frac{2}{3}$, which has not been improved since, in spite of many efforts, \cite{Zust}, \cite{LDZ}, \cite{Zust2016}. 

More generally, let $X$ be a metric space which is a topological manifold. Let $\alpha(X)$ be the supremum of $\alpha$ such that there exists a $C^{\alpha}$ homeomorphism from an open set in Euclidean space $\R^{\mathrm{dim}(X)}$ to $X$.

The case of Carnot groups equipped with Carnot-Carath\'eodory metrics is especially interesting, since they are prototypes of isometry-transitive and self-similar geodesic metric spaces. Let $\mathfrak{g}$ be a Carnot Lie algebra, i.e. a graded Lie algebra $\mathfrak{g}=\mathfrak{g}_1\oplus\cdots \mathfrak{g}_r$ with $[\mathfrak{g}_i,\mathfrak{g}_j]\subset \mathfrak{g}_{i+j}$, and such that $\mathfrak{g}_1$ generates $\mathfrak{g}$. Let $G$ denote the corresponding simply connected Lie group. A choice of Euclidean norm on $\mathfrak{g}_1$ determines a left-invariant Carnot-Carath\'eodory metric on $G$, of Hausdorff dimension
\begin{eqnarray*}
Q=\sum_{i=1}^{r}i\,\mathrm{dim}(\mathfrak{g}_i).
\end{eqnarray*}
Any two such metrics are bi-Lipschitz equivalent, therefore the following exponent is well-defined: let $\alpha(G)$ be the supremum of $\alpha$ such that there exists a $C^{\alpha}$ homeomorphism from an open set in Euclidean space $\R^{\mathrm{dim}(G)}$ to $G$. The exponential map is $\frac{1}{r}$-H\"older continuous, showing that $\alpha(G)\geq \frac{1}{r}$. 

Gromov has developed numerous tools to get lower bounds on $\alpha$, and has obtained the following results, among others.

\begin{stheorem*}[Gromov 1996]
Let $X$ be a metric space which is a topological manifold of dimension $n$ and Hausdorff dimension $Q$. Then $\alpha(X)\leq\frac{n}{Q}$.

Let $G$ be a Carnot group of topological dimension $n$ and Hausdorff dimension $Q$. Then $\alpha(G)\leq\frac{n-1}{Q-1}$.

Let $G=Heis_{\C}^m$ be the $m$-th Heisenberg group, of dimension $2m+1$. Then $\alpha(Heis_{\C}^m)\leq\frac{m+1}{m+2}$.
\end{stheorem*}

\subsection{The quasisymmetric H\"older equivalence problem}

We address a variant of the H\"older equivalence problem motivated by Riemannian geometry.

Let $M$ be a Riemannian manifold. Let $-1\leq\delta<0$. Say $M$ is {\em $\delta$-pinched} if sectional curvature ranges between $-a$ and $\delta a$ for some $a>0$. 
Define the {\em optimal pinching} $\delta(M)$ of $M$ as the least $\delta\geq-1$ such that $M$ is quasiisometric to a $\delta$-pinched complete simply connected Riemannian manifold.  

There are a few homogeneous Riemannian manifolds whose optimal pinching is known. \cite{pinch} deals with semi-direct products $\R\ltimes\R^{n}$ where $\R$ acts on $\R^n$ by matrices with only two distinct eigenvalues.

Optimal pinching is not known for the oldest examples, rank one symmetric spaces of noncompact type. These are hyperbolic spaces over the reals $H^{n}_{\R}$, the complex numbers $H^{m}_{\C}$, the quaternions $H^{m}_{\H}$, and the octonions $H^{2}_{\O}$. Real hyperbolic spaces have sectional curvature $-1$, and therefore optimal pinching $-1$. All other rank one symmetric spaces are $-\frac{1}{4}$-pinched. The optimal pinching of $H^{m}_{\C}$, $H^{m}_{\H}$ or $H^{m}_{\O}$ is conjectured to be $-\frac{1}{4}$, but still unknown. \cite{plan} suggests a different, yet nonconclusive, approach in the case of $H^{2}_{\C}$.

Negatively curved manifolds $M$ come with an ideal boundary $\partial M$, equipped with a family of (pairwise equivalent) \emph{visual metrics}. If curvature is $\leq -1$, visual metrics are true distances, i.e. satisfy triangle inequality. Polar coordinates are defined globally, extending into a homeomorphism from the round sphere to the ideal boundary, $\exp:S\to \partial M$. Rauch's comparison theorem shows that the behaviour of visual metrics in such coordinates reflects the distribution of sectional curvature: if $M$ is $\delta$-pinched, $\exp$ is distance increasing and $C^{\alpha}$-H\"older continuous with $\alpha=\sqrt{-\delta}$.

For instance, ideal boundaries of rank one symmetric spaces are round spheres or subRiemannian manifolds quasisymmetrically equivalent to Heisenberg groups $Heis_{\C}^{m-1}$, $Heis_{\H}^{m-1}$, $Heis_{\O}^{1}$. 

Finally, a quasiisometry between negatively curved manifolds $M$ and $M'$ extends to a quasisymmetric homeomorphism between ideal boundaries $\partial M\to\partial M'$.

Piling things up, a quasiisometry of a rank one symmetric space $M$ with a $\delta$-pinched manifold $M'$ provides us with a metric space $X'=\partial M'$ which is quasisymmetrically equivalent to a sub-Riemannian manifold and a $C^{\sqrt{-\delta}}$ homeomorphism from the round sphere $\exp:S\to X'$. Furthermore, $\exp^{-1}$ is Lipschitz.

\begin{defi}
Let $X$ be a metric space which is a topological manifold. Let $\alpha_{qs}(X)$ be the supremum of $\alpha$ such that there exists a metric space $X'$ quasisymmetric to an open subset of $X$ and a homeomorphism of an open subset of Euclidean space to $X'$ which is $C^\alpha$ H\"older continuous, and whose inverse is Lipschitz.
\end{defi}

\begin{que}[Quasisymmetric H\"older-Lipschitz equivalence problem]

Let $G$ be a non abelian Carnot group. Prove that $\alpha_{qs}(G)\leq\frac{1}{2}$.
\end{que}

The fact that Carnot groups of Haudorff dimension $Q$ have conformal dimension $Q$ (see \cite{confdim}) implies that $\alpha_{qs}(G)\leq \frac{n}{Q}$. We make one step further.

\begin{thm}
\label{thmintro}
Let $G$ be a Carnot group of topological dimension $n$ and Hausdorff dimension $Q$. Then $\alpha_{qs}(G)\leq\frac{n-1}{Q-1}$.
\end{thm}

This has an unsharp consequence for the pinching problem.

\begin{cor}
If $H^{2}_{\C}$ is quasiisometric to a $\delta$-pinched Riemannian manifold, then $\delta\geq -\frac{4}{9}$. 
\end{cor}
There are similar conclusions for other rank one symmetric spaces.

\section{The method}

Following Gromov, \cite{Gromovcarnot}, we manage to produce a smooth hypersurface in $\R^n$ whose image in the unknown space $X'$ has Hausdorff dimension at least $Q-1$. If $X'$ is a Carnot group, one knows, thanks to the isoperimetric inequality, that every hypersurface has Hausdorff dimension at least $Q-1$. The isoperimetric inequality or the Hausdorff dimension of subsets are not quasisymmetry invariants, so a different argument is needed.

Ours is based on a consequence of the coarea formula. If $X$ is a connected open subset of a Carnot group of Hausdorff dimension $Q$, and $u:X\to\R$ is a Lipschitz function, then (see \cite{Magnani}),
\begin{eqnarray}
\label{coareaineq}
\int_{X}|\nabla u|^Q\leq\int_{\R}\left(\int_{u^{-1}(t)}|\nabla u|^{Q-1}\right)\,dt\leq \mathrm{const.}\,\int_{\R}\mathcal{H}^{Q-1}(u^{-1}(t))\,dt.
\end{eqnarray}
Here, $\nabla u$ denotes the horizontal gradient. Since, for non constant $u$, $\int_{X}|\nabla u|^Q>0$, this shows that there exists $t\in\R$ such that $\mathcal{H}^{Q-1}(u^{-1}(t))>0$, and therefore $u^{-1}(t)$ has Hausdorff dimension at least $Q-1$. 

On the unknown metric space $X'$, we define, mimicking Tricot's construction of packing measures, quasisymmetry invariant analogues of the conformally invariant integrals $\int_{X}|\nabla u|^Q$ and $\int_{u^{-1}(t)}|\nabla u|^{Q-1}$. The sequence of inequalities (\ref{coareaineq}) persists (packing measures instead of covering measures are required in order to have the inequality in the appropriate direction). The main point is non vanishing of the packing analogue of $\int_{X}|\nabla u|^Q$. This is proved by a coarse analogue of the length area method (the covering version of which dates back to the 1980's, \cite{confdim}).

There is some hope to improve exponents and reach the sharp exponent $\frac{1}{2}$, at least for the simpler H\"older-Lipschitz equivalence problem, see the last section.

\textbf{Acknowledgements}. Many thanks to Zoltan Balogh, Artem Kozhevnikov, Enrico Le Donne and Roger Z\"ust for innumerable discussions about the H\"older equivalence problem.

\section{Packings}

\subsection{Bounded multiplicity packings}

Let $X$ be a metric space. A ball $B$ in $X$ is the data of a point $x\in X$ and a radius $r\geq 0$. For brevity, we also denote $B(x,r)$ by $B$. If $\lambda\geq 0$, $\lambda B$ denotes $B(x,\lambda r)$. 

\begin{defi}
Let $N$ be an integer, let $\ell\geq 1$. Let $X$ be a metric space. An \emph{$\ell$-packing} is a countable collection of balls $\{B_j\}$ such that concentric balls $\ell B_j$ are pairwise disjoint.
An \emph{$(N,\ell)$-packing} is a collection of balls $\{B_j\}$ which is the union of at most $N$ $\ell$-packings.
\end{defi}

\subsection{Covering and packing measures}

\begin{defi}
\label{defphi}
Let $\phi$ be a positive function on the set of balls in $X$. Two types of (pre-)measures can be obtained from it as follows. Let $A$ be a subset of $X$, $\epsilon>0$ and $p\geq 1$.

1. Covering measure
\begin{eqnarray*}
K\Phi_{N,\ell}^{p;\epsilon}(A)=\inf\{\sum_{i}\phi(B_{i})^p\,&;&\,\{B_{i}\}\,(N,\ell)-\mathrm{packing~of }\,X\,\mathrm{that~covers}\,A,\\
&&B_{i}\,\mathrm{of~diameter}\,\leq\epsilon\},
\end{eqnarray*}
\begin{eqnarray*}
K\Phi_{N,\ell}^{p}(A)=\lim_{\epsilon\to 0}K\Phi_{N,\ell}^{p;\epsilon}(A).
\end{eqnarray*}

2. Packing measure
\begin{eqnarray*}
P\Phi_{N,\ell}^{p;\epsilon}(A)=\sup\{\sum_{i}\phi(B_{i})^p\,&;&\,\{B_{i}\}\,(N,\ell)-\mathrm{packing~of }\,X,\\
&&B_{i}\,\mathrm{centered~on}\,A,\,\mathrm{of~diameter}\,\leq\epsilon\},
\end{eqnarray*}
\begin{eqnarray*}
P\Phi_{N,\ell}^{p}(A)=\lim_{\epsilon\to 0}P\Phi_{N,\ell}^{p;\epsilon}(A).
\end{eqnarray*}
\end{defi}

\begin{rem}
$K\Phi_{N,\ell}^{p}$ is a measure. $P\Phi_{N,\ell}^{p}$ is merely a pre-measure. To turn it into a measure, it suffices to force $\sigma$-additivity by declaring
\begin{eqnarray*}
\tilde{P}\Phi_{N,\ell}^{p}(A)=\inf\{\sum_{j\in\N}P\Phi_{N,\ell}^{p}(A_j)\,;\,A\subset \bigcup_{j\in\N}A_j\}.
\end{eqnarray*}
We shall ignore this point, we do not need $P\Phi_{N,\ell}^{p}$ to be a measure.
\end{rem}

\begin{exa}
\label{defpackmeas}
Take
\begin{eqnarray*}
\phi(B)=\delta(B)=\mathrm{radius}(B).
\end{eqnarray*}
The resulting covering measure $K\Delta_{N,\ell}^{p}$ is a minor variant of Hausdorff spherical measure. The packing pre-measure $P\Delta_{N,\ell}^{p}$ is called here the \emph{$p$-dimensional packing pre-measure} with parameters $(N,\ell)$.
\end{exa}

\begin{defi}
Let $N\in\N$, $\ell\geq 1$. Let $X$ be a metric space. The \emph{packing dimension} of a subset $A\subset X$ is
\begin{eqnarray*}
\mathrm{dim}_{N,\ell}(A)=\sup\{p\,;\,P\Delta_{N,\ell}^{p}(A)>0\}.
\end{eqnarray*}
\end{defi}

\begin{exa}
\label{segment}
Let $X$ be Euclidean space or a Carnot group. Let $A$ be a horizontal line segment of length $L$. For all $N\geq 1$ and $\ell\geq 1$, $P\Delta_{N,\ell}^{p;\epsilon}(A)\sim L\epsilon^{p-1}$. It follows that $\mathrm{dim}_{N,\ell}(A)=1$.
\end{exa}
Indeed, for any $(N,\ell)$-packing centered on $A$, $\sum \mathrm{radius}(B_i)\leq \frac{1}{2}NL$, so, if balls have radii $\leq \epsilon$, 
$$\sum \mathrm{radius}(B_i)^p\leq \frac{1}{2}NL\epsilon^{p-1}.$$
Conversely, there is an $(N,1)$-packing by $\epsilon$-balls which achieves this bound.

\begin{exa}
\label{vertical}
Let $X$ be Heisenberg group. Let $A$ be a vertical line segment of height $h$. For all $N\geq 1$ and $\ell\geq 1$, $P\Delta_{N,\ell}^{p;\epsilon}(A)\sim h\epsilon^{p-2}$. It follows that $\mathrm{dim}_{N,\ell}(A)=2$.
\end{exa}
Indeed, for any $(N,\ell)$-packing centered on $A$, $\sum \mathrm{radius}(B_i)^2\leq cNh$, where $c$ is the constant such that the height of an $r$-ball equals $\sqrt{\frac{r}{c}}$, so, if balls have radii $\leq \epsilon$, 
$$\sum \mathrm{radius}(B_i)^p\leq cNh\epsilon^{p-2}.$$
Conversely, there is an $(N,1)$-packing by $\epsilon$-balls which achieves this bound.

\begin{exa}
\label{hyperplane}
Let $X$ and $X'$ be open subsets of Carnot groups $G$ and $G'$ of Hausdorff dimensions $Q$ and $Q'$. Let $u:G\to G'$ be a surjective homogeneous homomorphism, with kernel $V$. Let $\mathcal{L}$ denote the Lebesgue measure on $V$. Let $A$ be a bounded open subset of $V$. Assume that the relative boundary $\partial A$ of $A$ in $V$ has vanishing Lebesgue measure $\mathcal{L}(\partial A)=0$. For all $\ell\geq 1$ and large enough $N$, $P\Delta_{N,\ell}^{p;\epsilon}(A)\sim \mathcal{L}(A)\epsilon^{p-Q+Q'}$. It follows that $\mathrm{dim}_{N,\ell}(A)=Q-Q'$.
\end{exa}
Indeed, by translation and dilation invariance, if $B$ is an $r$-ball centered on $A$, the Lebesgue measure $\mathcal{L}(A\cap B)=\frac{1}{c}\,r^{Q-Q'}$ for some constant $c>0$. For any $(N,\ell)$-packing by balls centered on $A$, of radius $\leq \epsilon$, 
$$\sum \mathrm{radius}(B_i)^{Q-Q'}\leq cN\mathcal{L}(A+\epsilon),$$
where $A+\epsilon$ denotes the $\epsilon$-neighborhood of $A$ in $V$. Thus
$$\sum \mathrm{radius}(B_i)^p\leq cN\mathcal{L}(A+\epsilon)\epsilon^{p-Q+Q'}.$$
Conversely, for large enough $N$, there is an $(N,1)$-packing by $\epsilon$-balls which achieves asymptotically this bound.

\subsection{Quasisymmetric invariance}
\label{qs}

\begin{prop}
\label{phiqs}
Let $X$, $X'$ be metric spaces. Let $f:X\to X'$ be a quasisymmetric homeomorphism. Given a function $\psi$ on balls of $X'$, define a function $\phi$ on balls of $X$ as follows. Given a ball $B$ of $X$ centered at $x$, let $B'$ be the smallest ball centered at $f(x)$ containing $f(B)$, and let $\phi(B)=\psi(B')$. Then, for all $\ell'\geq 1$, there exists $\ell\geq 1$ such that for all $N\in\N$, all $p>0$, and all subsets $A\subset X$,
\begin{eqnarray*}
K\Psi_{N,\ell'}^{p}(f(A)) &\leq& K\Phi_{N,\ell}^{p}(A),\\
P\Phi_{N,\ell}^{p}(A)&\leq& P\Psi_{N,\ell'}^{p}(f(A)).
\end{eqnarray*}
Symmetric statements also hold since $f^{-1}$ is quasisymmetric as well.
\end{prop}

\begin{pf}
Let $\eta:(0,+\infty)\to(0,+\infty)$ be the homeomorphism measuring the quasisymmetry of $f$, i.e. for every triple $x,y,z$ of distinct points of $X$,
\begin{eqnarray*}
\frac{d(f(x),f(y))}{d(f(x),f(z))}\leq \eta(\frac{d(x,y)}{d(x,z)}).
\end{eqnarray*}
Fix $\ell'\geq 1$ and $N\in\N$. Let $B$ be a ball of $X$. The procedure above defines a corresponding $B'$ of radius $\rho'$. Let $y\in B$ be such that $d(f(x),f(y))=\rho'$. If $z\in f^{-1}(\ell' B')$, $d(f(x),f(z))\leq \ell'\rho'$, so
\begin{eqnarray*}
\frac{d(f(x),f(y))}{d(f(x),f(z))}\geq \frac{\rho'}{\ell'\rho'}=\frac{1}{\ell'}.
\end{eqnarray*}
By quasi-symmetry, this implies that $\eta(\frac{d(x,y)}{d(x,z)})\geq \frac{1}{\ell'}$, and thus 
$$d(x,z)\leq \frac{1}{\eta^{-1}(\frac{1}{\ell'})}d(x,y).$$
In other words, $z\in \ell B$ with
$\ell=\frac{1}{\eta^{-1}(\frac{1}{\ell'})}$. We conclude that, for every $\lambda'\geq 1$, there exist $\lambda\geq 1$ such that $f(B)\subset B'$ and $f^{-1}(\lambda'B')\subset\lambda B$. We apply it with $\lambda'=3\ell'$ to get $\ell=\lambda\geq 1$.

We show that $(N,\ell)$-packings are mapped to $(N,\ell')$-packings by contradiction. Let $\{B_j\}$ be an $(N,\ell)$-packing of $X$. Let $B'_1,\ldots,B'_N$ be distinct balls from the corresponding packing of $X'$, ordered by increasing radii. Assume that $\bigcap_{j=1}^{N}\ell'B'_j$ is non empty. Since $\ell'B'_1$ and $\ell'B'_j$ intersect and $\ell'B'_1$ is smaller, $\ell'B'_1 \subset 3\ell'B'_j$, so $f(x_1)$ belongs to $\bigcap_{j=1}^{N}3\ell'B'_j$. Then
\begin{eqnarray*}
x_1 \in f^{-1}(\bigcap_{j=1}^{N}3\ell'B'_j)\subset \bigcap_{j=1}^{N}\ell B_j =\emptyset,
\end{eqnarray*}
contradiction. We conclude that $\{B'_j\}$ is an $(N,\ell')$-packing of $X'$, for $\ell=\frac{1}{\eta^{-1}(\frac{1}{3\ell'})}$. If $\{B_i\}$ covers a subset $A$, then $\{B'_i\}$ covers $f(A)$. If $\{B_i\}$ is centered on $A$, then $\{B'_i\}$ is centered on $f(A)$. So the announced inequalities follow from taking an infimum or a supremum.
\end{pf}

\subsection{H\"older covariance}

\begin{prop}
\label{holder}
Let $X$, $Y$ be metric spaces. Let $0<\alpha\leq 1$. Let $N\in\N$ and $\ell\geq 1$. Let $f:X\to Y$ be $C^\alpha$-H\"older continuous. Then, for every subset $A\subset X$,
\begin{eqnarray*}
P\Delta_{N,\ell}^{\alpha p}(A)\geq \mathrm{const.}\,P\Delta_{N,\ell}^{p}(f(A)).
\end{eqnarray*}
It follows that $\mathrm{dim}_{N,\ell}(A)\geq\alpha\,\mathrm{dim}_{N,\ell}(f(A))$.
\end{prop}

\begin{pf}
Let $\{B'_i\}$ be an $(N,\ell)$-packing of $Y$ by balls of radii $\leq\epsilon$ centered on $f(A)$. By assumption, for all $x$, $x'\in X$,
\begin{eqnarray*}
d(f(x),f(x'))\leq C\,d(x,x')^\alpha.
\end{eqnarray*}
If $B'_i=B(x'_i,r'_i)$, pick an inverse image $x_i\in f^{-1}(x'_i)\cap A$ and set 
\begin{eqnarray*}
B_i=B(x_i,\frac{1}{\ell}(\frac{\ell r'_i}{C})^{1/\alpha}).
\end{eqnarray*}
Then $B_i$ are balls centered on $A$, of radii $\leq\epsilon':=\frac{1}{\ell}(\frac{\ell \epsilon}{C})^{1/\alpha}$. Furthermore, the collection of sets
\begin{eqnarray*}
f(\ell B_i)\subset\ell B'_i
\end{eqnarray*}
has multiplicity $<N$, so the collection of balls $\{B_i\}$ is an $(N,\ell)$-packing of $X$. By definition,
\begin{eqnarray*}
P\Delta_{N,\ell}^{\alpha p;\epsilon'}(A) &\geq&\sum_{i}\mathrm{radius}(B_i)^{\alpha p}\\
&=&\sum_{i}(\frac{1}{\ell}(\frac{\ell r'_i}{C})^{1/\alpha})^{\alpha p}\\
&=&\ell^{p(1-\alpha)}C^{-p}\sum_{i}{r'}_i^{p}.
\end{eqnarray*}
Taking the supremum over $(N,\ell)$-packings yields
\begin{eqnarray*}
P\Delta_{N,\ell}^{\alpha p;\epsilon'}(A) \geq \ell^{p(1-\alpha)}C^{-p}P\Delta_{N,\ell}^{p;\epsilon}(f(A)). 
\end{eqnarray*}
Letting $\epsilon$ tend to 0, one concludes
\begin{eqnarray*}
P\Delta_{N,\ell}^{\alpha p}(A) \geq \ell^{p(1-\alpha)}C^{-p}P\Delta_{N,\ell}^{p}(f(A)). 
\end{eqnarray*}
\end{pf}

\section{The coarea inequality}

\subsection{Energy}

\begin{defi}
\label{defeu}
Let $X$ be a metric space. Let $u:X\to M$ be a map to an auxiliary measure space $(M,\mu)$. Let 
\begin{eqnarray*}
e_{u}(B)=\mu(u(B)).
\end{eqnarray*}
The total mass $PE_{u,N,\ell}^{p}(X)$ of the resulting packing pre-measure is the \emph{$p$-energy} of $u$ with parameters $(N,\ell)$.
\end{defi}

This is again a quasisymmetry invariant. If $f:X\to X'$ is a quasisymmetric homeomorphism, a map $u:X\to M$ gives rise to a map $u'=u\circ f^{-1}:X'\to M$. As in subsection \ref{qs}, given a ball $B$ of $X$ centered at $x$, let $B'$ be the smallest ball centered at $f(x)$ containing $f(B)$. Then $B\subset f^{-1}(B')$, $e_{u'}(B')=\mu(u(f^{-1}(B')))\geq \mu(u(B))=e_u(B)$, so we can assert that, for all $\ell'\geq 1$, there exists $\ell\geq 1$ such that for all $N\in\N$, all $p>0$ and all subsets $A\subset X$,
\begin{eqnarray*}
PE_{u,N,\ell}^{p}(A)\leq PE_{u',N,\ell'}^{p}(f(A)).
\end{eqnarray*}
Conversely, for all $\ell\geq 1$, there exists $\ell'\geq 1$ such that for all $N\in\N$, all $p>0$ and all subsets $A\subset X$,
\begin{eqnarray*}
PE_{u',N,\ell'}^{p}(f(A))\leq PE_{u,N,\ell}^{p}(A).
\end{eqnarray*}

\subsection{Coarea inequality}

\begin{prop}
\label{coarea}
Let $X$ be a metric space. Let $u:X\to M$ be a map to a measure space $(M,\mu)$. 
Let $N\in\N$ and $\ell\geq 1$. Then
\begin{eqnarray*}
PE_{u,N,2\ell}^{p}(X)\leq \int_{M}PE_{u,N,\ell}^{p-1}(u^{-1}(m))\,d\mu(m).
\end{eqnarray*}
\end{prop}

\begin{pf}
Let $\{B_i\}$ be an $(N,2\ell)$-packing of $X$ consisting of balls with radius $\leq \epsilon$. Write
\begin{eqnarray*}
\mu(u(B_i))=\int_{M}1_{u(B_i)}(m)\,dm.
\end{eqnarray*}
Furthermore, when $m\in u(B_i)$, pick a point $x_i\in B_i\cap u^{-1}(m)$ and let $B_{i,m}$ be the smallest ball centered at $x_i$ which contains $B_i$. Note that $B_{i,m}\subset 2B_i$, so that, for each $m\in M$, the collection $\{B_{i,m}\}$ is an $(N,\ell)$-packing of $X$ consisting of balls with radius $\leq 2\epsilon$ centered on $u^{-1}(m)$.
\begin{eqnarray*}
\sum_{i}\mu(u(B_i))^p 
&=&\sum_{i}(\int_{M}1_{u(B_{i})}(m)\,d\mu(m))\mu(u(B_i))^{p-1}\\
&=&\int_{M}(\sum_{i}1_{u(B_{i})}(m)\mu(u(B_i))^{p-1})\,d\mu(m)\\
&=&\int_{M}(\sum_{\{i\,;\,m\in u(B_{i})\}}\mu(u(B_i))^{p-1})\,d\mu(m)\\
&\leq&\int_{M}(\sum_{\{i\,;\,m\in u(B_{i})\}}\mu(u(B_{i,m}))^{p-1})\,d\mu(m)\\
&\leq&\int_{M}PE_{u,N,\ell}^{p-1;2\epsilon}(u^{-1}(m))\,d\mu(m).
\end{eqnarray*}
Taking the supremum over $(N,2\ell)$-packings yields
\begin{eqnarray*}
PE_{u,N,2\ell}^{p;\epsilon}(X) \leq \int_{M}PE_{u,N,\ell}^{p-1;2\epsilon}(u^{-1}(m))\,d\mu(m). 
\end{eqnarray*}
Letting $\epsilon$ tend to 0, one concludes that
\begin{eqnarray*}
PE_{u,N,2\ell}^{p}(X) \leq \int_{M}PE_{u,N,\ell}^{p-1}(u^{-1}(m))\,d\mu(m). 
\end{eqnarray*}
\end{pf}

\begin{rem}
Covering measures satisfy the opposite inequality
\begin{eqnarray*}
KE_{u,N,\ell}^{p}(X)\geq \int_{M}KE_{u,N,\ell}^{p-1}(u^{-1}(m))\,d\mu(m).
\end{eqnarray*}
\end{rem}

\section{Lower bounds on energy}

\subsection{Modulus estimate}

This is a packing version of a classical coarse modulus estimate for covering measures.

\begin{prop}
\label{modulus}
Let $X$ be a metric space. Let $\Gamma$ be a family of subsets of $X$, equipped with a measure $d\gamma$. For each $\gamma\in\Gamma$, a probability measure $m_{\gamma}$ is given on $\gamma$. Let $p\geq 1$, $N\in\N$, $\ell\geq 1$. Assume that there exist arbitrarily fine $(N,\ell)$-packings of $X$ that cover $X$. Assume that there exists a constant $\tau$ such that for every small enough ball $B$ of $X$,
\begin{eqnarray*}
\int_{\{\gamma\in\Gamma\,;\,\gamma\cap B\not=\emptyset\}}m_{\gamma}(\gamma\cap \ell B)^{1-p}\,d\gamma\leq \tau.
\end{eqnarray*}
Then, for every function $\phi$ on the set of balls of $X$,
\begin{eqnarray*}
P\Phi_{N,\ell}^{p}(X)\geq\frac{1}{N^{p-1}\tau}\int_{\Gamma}K\Phi_{N,\ell}^{1}(\gamma)^p\,d\gamma.
\end{eqnarray*}
\end{prop}

\begin{pf}
Let $\{B_i\}$ be an $(N,\ell)$-packing of $X$ by balls of radius $\leq\epsilon$ that covers $X$. 
Let $1_i$ be the function defined on $\Gamma$ by
\begin{eqnarray*}
1_i(\gamma)=\begin{cases}
1& \text{if }\gamma\cap B_i\not=\emptyset, \\
0& \text{otherwise}.
\end{cases}
\end{eqnarray*}
For each set $\gamma$, the balls such that $1_i(\gamma)=1$ cover $\gamma$, thus
\begin{eqnarray*}
K\Phi_{N,\ell}^{1;\epsilon}(\gamma)&\leq&\sum_{i}1_{i}(\gamma)\phi(B_i)\\
&=&\sum_{i}1_{i}(\gamma)\phi(B_i)m_\gamma(\gamma\cap\ell B_i)^{\frac{1-p}{p}}m_\gamma(\gamma\cap\ell B_i)^{\frac{p-1}{p}}.
\end{eqnarray*}
H\"older's inequality gives
\begin{eqnarray*}
K\Phi_{N,\ell}^{1;\epsilon}(\gamma)^p&\leq&\left(\sum_{i}1_{i}(\gamma)\phi(B_i)^p m_\gamma(\gamma\cap\ell B_i)^{1-p}\right)\left(\sum_{i}m_\gamma(\gamma\cap\ell B_i)\right)^{p-1}.
\end{eqnarray*}
Since the covering $\{\ell B_i\}$ has multiplicity $<N$ and $m_\gamma$ is a probability measure, the rightmost factor is $<N^{p-1}$. Integrating over $\Gamma$ gives
\begin{eqnarray*}
\int_{\Gamma}K\Phi_{N,\ell}^{1;\epsilon}(\gamma)^p \,d\gamma&\leq&N^{p-1}\sum_{i}\phi(B_i)^p \left(\int_{\Gamma}1_{i}(\gamma)m_\gamma(\gamma\cap\ell B_i)^{1-p}\,d\gamma\right)\\
&\leq&N^{p-1}\tau\sum_{i}\phi(B_i)^p \\
&\leq&N^{p-1}\tau P\Phi_{N,\ell}^{p;\epsilon}(X).
\end{eqnarray*}
\end{pf}

\begin{rem}
\label{doubling}
If $X$ is doubling at small scales, fine covering $(N,\ell)$-packings exist with $N$ depending only on $\ell$.
\end{rem}
Indeed, pick a maximal packing by disjoint $\frac{\epsilon}{2}$ balls. Then the doubled balls cover. If two $\ell$ times larger balls $B(x,\ell\epsilon)$ and $B(x',\ell\epsilon)$ overlap, then $B(x',\frac{\epsilon}{2})\subset B(x,(2\ell+1)\epsilon)$. The number of such balls is bounded above by the number of disjoint $\frac{\epsilon}{2}$-balls in a $(2\ell+1)\epsilon$-ball, which, in a doubling metric space, is bounded above in terms of $\ell$ only. Of course, we need this doubling property only for $\epsilon$ small.

\begin{exa}
Let $X$ be a Carnot group of Hausdorff dimension $Q$. Let $\Gamma$ be a family of parallel horizontal unit line segments. Then, for all $\ell>1$ and suitable $N$, the assumptions of Proposition \ref{modulus} are satisfied with $p=Q$, $m_\gamma$ the length measure, $d\gamma$ the Lebesgue measure on a codimension 1 subgroup.
\end{exa}
Indeed, apply translation and dilation invariance.

\begin{defi}
Let $\Gamma$ be a family of subsets of $X$. Say a function $\phi$ on the set of balls of $X$ is {\em $\Gamma$-admissible} if for every $\gamma\in\Gamma$,
\begin{eqnarray*}
K\Phi_{N,\ell}^{1}(\gamma)\geq 1.
\end{eqnarray*}
The \emph{$(p,N,\ell)$-modulus} of the family $\Gamma$ is the infimum of $P\Phi_{N,\ell}^{p}(X)$ over all $\Gamma$-admissible functions $\phi$.
\end{defi}
Thus Proposition \ref{modulus} states a sufficient condition for a family of subsets to have positive $p$-modulus. $p$-modulus is quasisymmetry invariant in the following sense. Let $f:X\to X'$ be a quasisymmetric homeomorphism. Then for all $\ell'\geq 1$, there exists $\ell\geq 1$ such that 
\begin{eqnarray*}
M_{N,\ell}^{p}(\Gamma)\leq M_{N,\ell'}^{p}(f(\Gamma)),
\end{eqnarray*}
and for all $\ell\geq 1$, there exists $\ell'\geq 1$ such that 
\begin{eqnarray*}
M_{N,\ell'}^{p}(f(\Gamma))\leq M_{N,\ell'}^{p}(\Gamma).
\end{eqnarray*}

\begin{exa}
Let $X$ be a Carnot group of Hausdorff dimension $Q$. For every $\ell\geq 1$, there exists $N\in\N$ such that families of parallel horizontal unit line segments have positive $(Q,N,\ell)$-modulus.
\end{exa}

\subsection{Quasiconformal submersions}

\begin{defi}
Let $X$ and $Y$ be metric spaces. Say a continuous map $u:X\to Y$ is a \emph{quasiconformal submersion} if there exists a homeomorphism $\eta:(0,+\infty)\to(0,+\infty)$ such that for every ball $B$ of $X$, there exists a ball $B'$ of $Y$ such that for all $\lambda\geq 1$,
\begin{eqnarray*}
B'\subset u(B)\subset u(\lambda B)\subset\eta(\lambda)B'.
\end{eqnarray*}
\end{defi}
This notion is quasisymmetry invariant, both on the domain and on the range.

\begin{exa}
\label{rqs}
If balls in $X$ are connected and $Y=\R$, any continuous map $X\to Y$ is a quasiconformal submersion.
\end{exa}
By quasisymmetry invariance, this fact generalizes to spaces quasisymmetrically equivalent to metric spaces whose balls are connected.

\begin{exa}
Let $X=G$ and $Y=G'$ be Carnot groups. Then any surjective homogeneous homomorphism $u:G\to G'$ is a quasiconformal submersion. More generally, any contact map between open sets of Carnot groups whose differential is continuous and surjective is locally a quasiconformal submersion.
\end{exa}

\begin{cor}
\label{subm}
Let $X$ be a connected open subset of a Carnot group of Hausdorff dimension $Q$. Let $M$ be a metric measure space. Assume that there exist constants $d$ and $\nu$ such that for all small enough balls in $M$,
\begin{eqnarray*}
\mu(B(m,r))\geq \nu\,r^d.
\end{eqnarray*}
Let $u:X\to M$ be a nonconstant quasiconformal submersion. Then $$PE_{u,N,\ell}^{Q/d}(X)>0.$$
\end{cor}

\begin{pf}
Let $\phi=e_u^{1/d}$. Let $\eta$ be the function measuring quasiconformality of $u$. For every ball $B$ in $X$, $u(B)$ contains a ball of radius $\geq \frac{\mathrm{diameter(u(B))}}{\eta(1)}$, thus
\begin{eqnarray*}
\phi(B)=\mu(u(B))^{1/d}\geq \frac{\nu^{1/d}}{\eta(1)}\mathrm{diameter(u(B))}.
\end{eqnarray*}
Let $\Gamma\subset X$ be a family of parallel horizontal line segments of equal lengths. If a collection $\{B_i\}$ of balls covers one such segment $\gamma$, 
\begin{eqnarray*}
\mathrm{diameter}(u(\gamma))\leq\sum_{i}\mathrm{diameter}(u(\gamma\cap B_i))\leq \frac{\eta(1)}{\nu^{1/d}}\sum_{i}\phi(B_i).
\end{eqnarray*}
This shows that 
\begin{eqnarray*}
\int_{\Gamma}K\Phi_{N,\ell}^{1}(\gamma)^p\,d\gamma\geq\frac{\nu^{1/d}}{\eta(1)}\int_{\Gamma}\mathrm{diameter}(u(\gamma))^p\,d\gamma.
\end{eqnarray*}
Assume by contradiction that $PE_{u,N,\ell}^{Q/d}(X)=P\Phi_{N,\ell}^{Q}(X)=0$. Proposition \ref{modulus} implies that $u$ is constant on every segment $\gamma\in\Gamma$. This proves that $u$ is constant on every horizontal segment, hence on every polygonal curve made of horizontal segments. Since such curves allow to travel from any point to any other point of $X$, $u$ is constant, contradiction.
\end{pf}

\section{On the quasi-symmetric H\"older-Lipschitz problem}

\subsection{Energy dimension}

\begin{defi}
Let $X$ be a metric space. Define its \emph{energy dimension} as the infimum of exponents $p$ such that there exist non constant continuous functions $u:X\to\R$ with finite $p$-energy $PE_{u,N,\ell}^{p}(X)$ for $\ell=1$ and all large enough $N$.
\end{defi}
This is a quasisymmetry invariant. By definition, this is less than packing dimension. Think of it as an avatar of \emph{conformal dimension} (a generic term for the infimum of all dimensions of metric spaces quasisymmetrically equivalent to $X$, see \cite{confdim}, \cite{MKT}).

\begin{exa}
\label{exendim}
The energy dimension of an open subset of a Carnot group is equal to its Hausdorff dimension.
\end{exa}
Indeed, since real valued functions are quasiconformal submersions, this follows from Corollary \ref{subm}.

\begin{lem}
\label{submpositive}
Let $X$ be a metric space of energy dimension $Q$. Let $M$ be a $d$-Ahlfors regular metric space (at small scales). If $u:X\to M$ is a quasiconformal submersion, then for all $p<Q$, $PE_{u,N,\ell}^{p/d}(X)>0$. 
\end{lem}

\begin{pf}
$d$-Ahlfors regular at small scales means that for all small enough balls in $M$,
\begin{eqnarray*}
\nu\,r^d\leq \mu(B(m,r))\leq \frac{1}{\nu}\,r^d.
\end{eqnarray*}
Under this assumption, there exists a constant $c$ such that for all small enough balls $B$ of $X$,
\begin{eqnarray*}
\frac{1}{c}\mathrm{diameter}(u(B))\leq e_u(B)^{1/d}\leq c\,\mathrm{diameter}(u(B)).
\end{eqnarray*}
Let $p<Q$. We prove by contradiction that $PE_{u,N,\ell}^{p/d}(X)>0$. Pick a point $m_0\in M$. Define a real valued function $v$ on $X$ by
\begin{eqnarray*}
v(x)=d(m_0,u(x)).
\end{eqnarray*}
For all small enough balls $B$ in $X$,
\begin{eqnarray*}
\mathrm{diameter}(v(B))\leq \mathrm{diameter}(u(B))\leq c\,e_u(B)^{1/d},
\end{eqnarray*}
thus
\begin{eqnarray*}
PE_{v,N,\ell}^{p}(X)\leq c^p\,PE_{u,N,\ell}^{p/d}=0.
\end{eqnarray*}
By definition of energy dimension, this implies that $v$ is constant. Since this holds for every $m_0\in M$, one finds that $u$ is constant, contradiction. One concludes that $PE_{u,N,\ell}^{p/d}(X)>0$ for all $p<Q$.
\end{pf}

\subsection{Proof of the main theorem}

\begin{lem}
\label{qqq}
Let $X\subset G$ and $M\subset G''$ be open subsets of Carnot groups of Hausdorff dimensions $Q$ and $Q''$. Let $X'$ be a metric space. Let $f:X\to X'$ be a $C^\alpha$-H\"older continuous homeomorphism. Assume that $f^{-1}:X'\to X$ is Lipschitz. Let $u:G\to G''$ be a surjective homomorphism mapping $X$ to $M$. Let $u'=u\circ f^{-1}:X'\to M$. Assume that for all $p<\frac{Q'}{Q''}$, all $\ell>1$ and $N$ large enough, $PE_{u',N,2\ell}^{p}(X')>0$ . Then 
\begin{eqnarray*}
\alpha\leq\frac{Q-Q''}{Q'-Q''}.
\end{eqnarray*}
\end{lem}

\begin{pf}
Fix $\ell>1$ and choose $N$ according to Example \ref{doubling}.
For all $m\in M$, $\mathrm{dim}_{N,\ell}(u^{-1}(m))=Q-Q''$ (Example \ref{hyperplane}). By assumption, for all $p<Q'/Q''$, $PE_{u',N,2\ell}^{p}(X')>0$. Proposition \ref{coarea} implies that there exists $m_p\in M$ such that $PE_{u',N,\ell}^{p-1}(u'^{-1}(m_p))>0$. 

If $f^{-1}$ is Lipschitz, so is $u'$, thus $e_{u'}(B)\leq \mathrm{const.}\,\mathrm{radius}(B)^{Q''}$, so $PE_{u',N,\ell}^{p-1}\leq\mathrm{const.}P\Delta_{N,\ell}^{Q''(p-1)}$. It follows that $P\Delta_{N,\ell}^{Q''(p-1)}(u'^{-1}(m_p))>0$, and 
$$\mathrm{dim}_{N,\ell}(u'^{-1}(m_p))\geq Q''(p-1).$$ 
Since $u'^{-1}(m_p)=f(u^{-1}(m_p))$, H\"older covariance (Proposition \ref{holder}) implies that $Q-Q''\geq\alpha Q''(p-1)$. Since this holds for all $p<Q'/Q''$, $Q-Q''\geq\alpha(Q'-Q'')$.
\end{pf}

\begin{thm}
\label{thmconf}
Let $X\subset G$ be a connected open subset of a Carnot group of Hausdorff dimension $Q$. Let $X'$ be a metric space of energy dimension $Q'$. Assume that $X'$ is quasisymmetric to a metric space whose balls are connected. Let $f:X\to X'$ be a $C^\alpha$-H\"older continuous homeomorphism. Assume that $f^{-1}:X'\to X$ is Lipschitz. Then 
\begin{eqnarray*}
\alpha\leq\frac{Q-1}{Q'-1}.
\end{eqnarray*}
\end{thm}

\begin{pf}
Let $u:X\to\R$ be the restriction to $X$ of a non constant group homomorphism $G\to\R$. By assumption, $u'=u\circ f^{-1}:X'\to\R$ is automatically a quasiconformal submersion (Example \ref{rqs}). According to Lemma \ref{submpositive}, for all $p<Q'$, $PE_{u',N,2\ell}^{p}(X')>0$. Lemma \ref{qqq} implies that $Q-1\geq\alpha(Q'-1)$.
\end{pf}

Theorem \ref{thmintro} is the special case where $G=\R^n$ and $X'$ is quasisymmetrically equivalent to an open subset of a Carnot group of Hausdorff dimension $Q$. By quasisymmetry invariance and Example \ref{exendim}, $X'$ has energy dimension $\geq Q$, so all assumptions of Theorem \ref{thmconf} are satisfied.

\section{Speculation}

\subsection{Energies and homeomorphisms}

Lemma \ref{qqq} applies successfully with $G''=\R$, since we have some information on energy dimensions of Carnot groups. To get closer to the conjectured estimate $\alpha_{qs}(Heis_{\C}^1)\leq\frac{1}{2}$, one needs understand the energies of maps of Carnot groups $G'$ to $G''=\R^2$. 2-energy is especially relevant, since nonvanishing of 2-energy would lead to the sharp bound $\frac{1}{2}$.

Let us define the ``$\R^2$-dimension'' of a metric space as the infimum of exponents $p$ such that there exist non constant continuous \emph{open} maps $u:X\to\R^2$ with finite $p$-energy. The ``$\R^2$-dimension'' of Heisenberg group $Heis_{\C}^1$ is $\leq \frac{4}{3}$, thus much smaller than 2. Indeed, the projection $Heis_{\C}^1 \to\R^2$, $(x,y,z)\mapsto (y,z)$ whose fibers are horizontal line segments has finite $\frac{4}{3}$-energy. 

Hence one must better exploit the special characters of maps involved.

\begin{que}
Let $g:Heis^1_{\C}\to\R^3$ be a (local) homeomorphism. Does there exist a nonzero linear map $u:\R^3\to\R^2$ such that $u\circ g:Heis^1_{\C}\to\R^2$ has positive 2-energy ?
\end{que}

\begin{prop}
\label{diff}
Let $X'$ be an open subset of Heisenberg group $Heis^1_{\C}$. Let $g:X'\to\R^3$ be almost everywhere differentiable. Assume that for all linear maps $u:\R^3\to\R^2$, $PE_{u\circ g,N,\ell}^{2}(X')=0$. Then the differential $Dg\in Hom(Heis^1_{\C},\R^3)$ has rank $\leq 1$ almost everywhere.
\end{prop}

This suggests the following strategy:
\begin{enumerate}
  \item Let $X$ be a metric space. Prove that the composition $g$ of a quasisymmetric homeomorphism $Heis^1_{\C}\to X$ and a Lipschitz map $X\to\R^3$ is almost everywhere differentiable and absolutely continuous on lines.
  \item Prove that if $Dg$ has rank $\leq 1$ almost everywhere, then $g$ cannot be a homeomorphism. 
\end{enumerate}

\subsection{Proof of Proposition \ref{diff}} 

\begin{lem}\label{PEJ}
Let $X\subset G$, $X'\subset G'$ be open subsets of Carnot groups of Hausdorff dimensions $Q$ and $Q'$. Let $\beta$ denotes the unit ball in $G'$. If $h:G'\to G$ is a homogeneous homomorphism, define
\begin{eqnarray*}
J(h)=\frac{vol(h(\beta))}{vol(\beta)^{Q/Q'}},
\end{eqnarray*}
if $h$ is surjective, $J(h)=0$ otherwise. Let $g:X'\to X$ be a map which is almost everywhere differentiable. Then, for all $N$ and $\ell\geq 1$, 
\begin{eqnarray*}
PE_{g,N,2\ell}^{Q'/Q}(X')\geq C(N,\ell)\,\int_{X'}J(D_y g)^{Q'/Q}\,dy.
\end{eqnarray*}
\end{lem}

\begin{pf}
Let us first prove that at points $y\in X'$ where $g$ is differentiable,
\begin{eqnarray*}
\lim_{r\to 0}\frac{vol(g(B(y,r)))}{vol(B(y,r))^{Q/Q'}}= J(D_y g). 
\end{eqnarray*}
Up to translating, one can assume that $y$ is the identity element. Then $B(y,r)=\delta_r(\beta)$, 
\begin{eqnarray*}
\frac{vol(g(B(y,r)))}{vol(B(y,r))^{Q/Q'}}=\frac{vol(\delta_{1/r}g\delta_r(\beta))}{vol(\beta)^{Q/Q'}}.
\end{eqnarray*}
By definition of differentiability (\cite{annals}), $\delta_{1/r}g\delta_r$ converge uniformly to $D_y g$, hence the indicatrix $1_{\delta_{1/r}g\delta_r(\beta)}$ converges pointwise to $1_{D_y g(\beta)}$ away from the null-set $D_y g(\partial\beta)$, dominated convergence applies, and volumes converge.

Let $\beta_j$ denote the $1/j$-ball in $G'$. Let $Y\subset X'$ be an open subset compactly contained in $X'$. For $j$ large enough, define a measurable function $\eta_j$ on $Y$ as follows.
\begin{eqnarray*}
\eta_j(y)=j\sup_{z\in \beta_j}d(D_yg(z),g(y)^{-1}g(yz)).
\end{eqnarray*}
This function measures the speed at which $g$ is approximated by its differential. By assumption, as $j$ tends to $\infty$, $\eta_j$ tends to 0 almost everywhere. According to Lusin's theorem, the convergence is uniform on a compact set $Z$ whose complement has arbitrarily small measure. 

It follows that, at points $y\in Z$,
\begin{eqnarray*}
\frac{vol(g(y\beta_j))}{vol(y\beta_j)^{Q/Q'}}\to J(D_y g)
\end{eqnarray*}
uniformly. Therefore
\begin{eqnarray*}
vol(g(y\beta_j))^{Q'/Q}\sim J(D_y g)^{Q'/Q}vol(y\beta_j)\quad \textrm{as}\quad j\to\infty
\end{eqnarray*}
uniformly as $y$ varies in $Z$. Up to replacing balls containing points of $Z$ with twice larger balls centered at points of $Z$, one gets the same conclusion, with a loss of a power of 2, for balls containing a point of $Z$. Pick an $(N,\ell)$-packing of $X'$ of mesh $\leq\epsilon$, covering a subset of almost full measure of $X'$, and therefore a large part $Z'$ of $Z$. Discard balls which do not intersect $Z$. For the remaining balls,
\begin{eqnarray*}
J(D_y g)^{Q'/Q}vol(B)\lesssim vol(g(B))^{Q'/Q},
\end{eqnarray*}
therefore
\begin{eqnarray*}
\int_{Z'}J_\epsilon(y)^{Q'/Q}\,dy\lesssim N\sum_{B}vol(g(B))^{Q'/Q}\leq N\,PE^{Q'/Q}_{N,\ell}(X'),
\end{eqnarray*}
where 
\begin{eqnarray*}
J_\epsilon(z)=\inf_{B(z,\epsilon)}J(Dg).
\end{eqnarray*}
As $\epsilon$ decreases to 0, $J_\epsilon$ increases and converges almost everywhere to $J(Dg)$, so, by monotone convergence,
\begin{eqnarray*}
\int_{Z'}J(D_y g)^{Q'/Q}\,dy\leq C(N,\ell)\, PE^{Q'/Q}_{N,\ell}(X').
\end{eqnarray*}
Finally, the measure of $X'\setminus Z'$ is arbitrarily small, thus $Z'$ can be replaced with $X'$.
\end{pf}

Let $g$ be a map from an open subset of Heisenberg group to $\R^3$ which is differentiable almost everywhere. Assume that $Dg$ has rank 2 on a set $Y$ of positive measure. At each point of $Y$, one of the three projections $u:\R^3\to\R^2$ to coordinate planes is surjective. Thus, for one of them, $u$, this holds for a subset of positive measure. At such points, $J(D(u\circ g))>0$, thus $PE^2_{u\circ g,N,\ell}(Y)>0$.

\bigskip

\noindent
Pierre Pansu 
\par\noindent Laboratoire de Math\'ematiques d'Orsay,
\par\noindent Univ. Paris-Sud, CNRS, Universit\'e
Paris-Saclay
\par\noindent 91405 Orsay, France.
\par\noindent
e-mail: pierre.pansu@math.u-psud.fr

\end{document}